\documentclass[12pt,oneside,reqno]{amsart}
\usepackage{txfonts}
\usepackage{bbm}
\usepackage{amsmath}
\usepackage{graphicx}
\usepackage{mathrsfs}
\usepackage{stmaryrd}
\usepackage{amsfonts}
\usepackage{enumerate,amsmath,amssymb,amsthm}

\oddsidemargin=0mm \topmargin=-12mm \numberwithin{equation}{section}

\newcommand{\be}{\begin{eqnarray}}
\newcommand{\ee}{\end{eqnarray}}
\newcommand{\ce}{\begin{eqnarray*}}
\newcommand{\de}{\end{eqnarray*}}
\newtheorem{theorem}{Theorem}[section]
\newtheorem{lemma}[theorem]{Lemma}
\newtheorem{remark}[theorem]{Remark}
\newtheorem{definition}[theorem]{Definition}
\newtheorem{proposition}[theorem]{Proposition}
\newtheorem{Examples}[theorem]{Example}
\newtheorem{corollary}[theorem]{Corollary}

\def\e{{\mathrm{e}}}
\def\eps{\varepsilon}

\def\p{\partial}

\def\[{{\Big[}}
\def\]{{\Big]}}
\def\<{{\langle}}
\def\>{{\rangle}}
\def\({{\Big(}}
\def\){{\Big)}}

\def\bx{{\mathbf{x}}}

\def\dif{{\mathord{{\rm d}}}}

\def\no{\nonumber}
\def\={&\!\!=\!\!&}
\def\bt{\begin{theorem}}
\def\et{\end{theorem}}
\def\bl{\begin{lemma}}
\def\el{\end{lemma}}
\def\br{\begin{remark}}
\def\er{\end{remark}}

\def\bd{\begin{definition}}
\def\ed{\end{definition}}
\def\bp{\begin{proposition}}
\def\ep{\end{proposition}}
\def\bc{\begin{corollary}}
\def\ec{\end{corollary}}
\def\bx{\begin{Examples}}
\def\ex{\end{Examples}}

\def\cB{{\mathcal B}}

\def\cL{{\mathcal L}}

\def\mD{{\mathbb D}}
\def\mE{{\mathbb E}}

\def\mH{{\mathbb H}}

\def\mN{{\mathbb N}}

\def\mP{{\mathbb P}}

\def\mR{{\mathbb R}}
\def\mS{{\mathbb S}}

\def\mW{{\mathbb W}}

\def\sB{{\mathscr B}}

\def\sF{{\mathscr F}}

\def\sU{{\mathscr U}}
\def\sV{{\mathscr V}}

\def\geq{\geqslant}
\def\leq{\leqslant}

\allowdisplaybreaks

\begin{document}

\bf\title{Strong Feller properties for degenerate SDEs with jumps}

\date{}
\author{Zhao Dong, Xuhui Peng, Yulin Song and Xicheng Zhang}

\address{Zhao Dong:
Academy of Mathematics and Systems Sciences, Chinese Academy of Sciences, Beijing, 100190, P.R.China\\
Email: dzhao@amt.ac.cn
 }

\address{Xuhui Peng:
Academy of Mathematics and Systems Sciences, Chinese Academy of Sciences, Beijing, 100190, P.R.China\\
Email: pengxuhui@amss.ac.cn
 }

\address{Yulin Song:
School of Mathematical Sciences, Beijing Normal University, Beijing, {\rm 100875}, P.R.China\\
Email: songyl@amss.ac.cn
 }
 
\address{Xicheng Zhang:
School of Mathematics and Statistics, Wuhan University,
Wuhan, Hubei 430072, P.R.China\\
Email: XichengZhang@gmail.com
 }

\dedicatory{}

\begin{abstract}
Under full H\"ormander's conditions, we prove the strong Feller property of the semigroup determined by an SDE driven by additive subordinate Brownian motion,
where the drift is allowed to be arbitrarily growth. For this, we extend
a criterion due to Malicet-Poly \cite{Ma-Po} and Bally-Caramellino \cite{Ba-Ca} about 
the convergence of the laws of  Wiener functionals in total variations.
Moreover, the example of a chain of coupled oscillators is verified.
\end{abstract}

\maketitle

\rm


\section{Introduction}

Let $\mW$ be the space of all continuous functions from $\mR_+:=[0,\infty)$ to $\mR^m$ vanishing at starting point $0$, 
which is endowed with the locally uniform convergence topology and the Wiener measure $\mu_\mW$ so that the coordinate process
$$
W_t(\omega)=\omega_t
$$
is a standard $m$-dimensional Brownian motion. Let $\mH\subset\mW$ be the Cameron-Martin space consisting of all
absolutely continuous functions with square integrable derivatives. The inner product in $\mH$ is denoted by 
$$
\<h_1,h_2\>_\mH:=\sum_{i=1}^m\int^\infty_0\dot h^i_1(s)\dot h^i_2(s)\dif s.
$$
The triple $(\mW,\mH,\mu_{\mW})$ is also called the classical Wiener space.

Let $D$ be the Malliavin derivative operator. For $k\in\mN$ and $p\geq 1$,  let $\mD^{k,p}$ be the associated Wiener-Sobolev space with the norm:
$$
\|F\|_{k,p}:=\|F\|_p+\|DF\|_p+\cdots+\|D^kF\|_p,
$$
where $\|\cdot\|_p$ is the usual $L^p$-norm. Let $X:\mW\to\mR^d$ be a smooth Wiener functional in $\cap_{k,p} \mD^{k,p}$. 
Let $\Sigma^X_{ij}:=\<DX^i, DX^j\>_\mH$ be the Malliavin covariance matrix.
The classical Malliavin calculus studies the problem that under what conditions on $X$, the law of $X$ has a smooth density with
respect to the Lebesgue measure. In particular, as his theory's application, Malliavin gave a probabilistic proof for the celebrated H\"ormander's hypoellipticity theorem (cf. \cite{Ma, Nu}).
Nowadays, the Malliavin calculus, as a kind of infinite dimensional analysis, has been extensively used in many fields such as heat kernel estimates, 
large deviation theory, financial mathematics, numerical calculations, and so on (cf. \cite{Bi0} \cite{Ku-St1} \cite{Fo}).

 On the other hand, in the studies of the ergodicity of stochastic dynamical systems, the notion of strong Feller property plays a crucial role (cf. \cite{Da-Za}),
which relates to the following problem: Let $\Lambda$ be a metric space and
$(X_\lambda)_{\lambda\in\Lambda}$  a random field. We want to seek conditions on $X_\lambda$ so that for any  $f\in \cB_b(\mR^d)$
(the space of bounded measurable functions),
$$
\lambda\mapsto\mE f(X_\lambda)\mbox{ is continuous.}
$$
In many cases, this is difficult to verify. As we learned, if $X_t(x)$ is the solution of an SDE, there are many ways to derive the strong Feller property to 
$P_tf(x):=\mE f(X_t(x))$. For examples, Bismut-Elworthy-Li's formula provides an explicit formula for $\nabla P_tf(x)$ (cf. \cite{El-Li}).
Moreover, F.Y. Wang's Hanarck inequality gives some quantitive estimate to $P_tf(x)$ for finite and infinite dimensional systems, which can  also be used to derive the 
strong Feller property (cf. \cite{Wa}).

In the framework of the Malliavin calculus, the above problem can be introduced as follows. The celebrated Bouleau-Hirsch's criterion says that if 
$X_\lambda\in \mD^{1,p}$ for some $p>1$ and the Malliavin covariance matrix $\Sigma^X_\lambda:=\Sigma^{X_\lambda}$ is invertible almost surely, then the law of $X_\lambda$
is absolutely continuous with respect to the Lebesgue measure (cf. \cite{Nu}). But we have no any information about the regularity of the density $\rho_\lambda$.
In order to obtain such information, one usually needs the stronger hypothesis $(\Sigma^X_\lambda)^{-1}\in \cap_{p\geq 1}L^p$. If this is true, and if we work with a diffusion
process, then the semigroup of the diffusion has a ``regularization effect''. Question: is it possible to emphasis a regularization effect under the weak hypothesis 
``$\det(\Sigma^X_\lambda)>0$ almost surely''? The answer is yes. In fact, Bogachev \cite[Corollary 9.6.12]{Bo} has already shown the following result:
Let $X_n$ and $X$ be $d$-dimensional random variables in $\mD^{1,p}$ so that $X_n\to X$ in $\mD^{1,p}$. If $p\geq d$ and for almost all $\omega$,
$$
\{D_hX(\omega),\ h\in\mH\}=\mR^d,
$$
then the laws of $X_n$ converge to the law of $X$ in total variations. Notice that $\det(\Sigma^X(\omega))>0$ implies the above condition. This can be seen as follows:
Suppose that $\{D_hX(\omega),\ h\in\mH\}\not=\mR^d$, then there is a non-zero vector $v=(v_1,\cdots, v_d)\in\mR^d$ such that
$$
\<D_hX(\omega),v\>_{\mR^d}=0\ \ \forall h\in\mH\Rightarrow \sum_iv_iDX^i(\omega)=0\Rightarrow\Sigma^X(\omega)v=0
\Rightarrow \det(\Sigma^X(\omega))=0.
$$
This criterion recently was reproven by Malicet and Poly in \cite[Corollary 2.2]{Ma-Po} by using another argument (see also Bally and Caramellino \cite[Corollary 2.16]{Ba-Ca})

The first aim of this work is to extend Bogachev's result as follows.
\bt\label{Th1}
Let $(X_\lambda)_{\lambda\in\Lambda}$ be a family of  $\mR^d$-valued Wiener functionals over $\mW$. Suppose that for some $p>1$,
\begin{enumerate}[{\bf(H1)}]
\item $X_\lambda \in\mD^{2,p}$ for each $\lambda\in\Lambda$, and $\lambda\mapsto\|X_\lambda\|_{2,p}$ is locally bounded.
\item $\lambda\mapsto X_\lambda$ is continuous in probability, i.e., 
for any $\eps>0$ and $\lambda_0\in\Lambda$,
$$
\lim_{\lambda\to\lambda_0}\mP(|X_\lambda-X_{\lambda_0}|\geq\eps)=0.
$$
\item For each $\lambda\in\Lambda$, the Malliavin covariance matrix $\Sigma^X_\lambda$ of $X_\lambda$ is invertible almost surely.
\end{enumerate}
Then the law of $X_\lambda$ in $\mR^d$ admits a density $\rho_\lambda(x)$ so that $\lambda\mapsto \rho_\lambda$ is continuous in $L^1(\mR^d)$.
\et
\br
Our proof is different from \cite{Bo, Ma-Po, Ba-Ca} and based on the Sobolev's compact embedding. 
Compared with \cite{Bo}, our result requires less integrability and continuity, while more differentiability is needed. This can be considered as the case 
that the differentiability index can compensate
the integrability index in infinite dimensional calculus.
\er
Our another aim of this work is to apply the above criterion to the SDE driven by degenerate jump noises. 
Let $\mS$ be the space of all c\`adl\`ag functions from $\mR_+$ to $\mR^m_+$ with
$\ell_0=0$ and each component being increasing and purely jumping. Suppose that $\mS$
is endowed with the Skorohod metric and the probability measure $\mu_\mS$ so that the coordinate process 
$$
S_t(\ell):=\ell_t=(\ell^1_t,\cdots,\ell^m_t)
$$
is an $m$-dimensional L\'evy process with Laplace transform
\begin{align}
\mE^{\mu_{\mS}}(\e^{-z\cdot S_t})=\exp\left\{\int_{\mR^m_+}(\e^{-z\cdot u}-1)\nu_S(\dif u)\right\}.\label{EW111}
\end{align}
Consider the following product probability space
$$
(\Omega,\sF,\mP):=\Big(\mW\times \mS, \sB(\mW)\times\sB(\mS), \mu_\mW\times\mu_{\mS}\Big).
$$
If we lift $W_t$ and $S_t$ to this probability space, then $W_t$ and $S_t$ are independent, and the subordinated Brownian motion
$$
W_{S_t}:=\Big(W^1_{S^1_t},\cdots, W^m_{S^m_t}\Big)
$$
is an $m$-dimensional L\'evy process. Below we assume
\begin{align}
\mP(\omega\in\Omega: \exists j=1,\cdots,m  \mbox{ and } \exists t>0 \mbox{ such that } S^j_t(\omega)=0)=0,\label{ER}
\end{align}
which means that $S_t$ is nondegenerate along each direction. 

Consider the following SDE driven by $W_{S_t}$:
\begin{align}
\dif X_t=b(X_t)\dif t+A\dif W_{S_t},\ \ X_0=x,\label{SDE}
\end{align}
where $b:\mR^d\to\mR^d$ is a smooth function, $A=(a_{ij})$ is a $d\times m$ constant matrix.
Let $H:\mR^d\to\mR^+$ be a $C^\infty$-function with $\lim_{|x|\to\infty}H(x)=\infty$, which is called a Lyapunov function.
We assume that for some Lyapunov function $H$ and $\kappa_1,\kappa_2,\kappa_3\geq0$,
\begin{align}
b(x)\cdot\nabla H (x)&\leq \kappa_1 H(x),\label{EE1}
\end{align}
and for all $k=1,\cdots,m$,
\begin{align}
\Big|\sum_i\p_i H(x)a_{ik}\Big|^2\leq \kappa_2H (x),\ \ \ \sum_{ij}\p_i\p_j H(x)a_{ik}a_{jk}\leq \kappa_3.\label{EE3}
\end{align}
Under (\ref{EE1})-(\ref{EE3}), X. Zhang in \cite[Theorem 3.1]{Zh3} has already proved that SDE (\ref{SDE}) has a unique solution $X_t(x)$, which defines a Markov process.
The associated Markov semigroup is defined by
$$
P_tf(x):=\mE f(X_t(x)).
$$
We say that $(b,A)$ satisfies a H\"ormander's condition at point $x\in\mR^d$ if for some $n=n(x)\in \mathbb{N}$,
\begin{eqnarray}
\text{Rank}[A,B_1(x)A,B_2(x)A,\cdots,B_n(x)A]=d,\label{eq}
\end{eqnarray}
where $B_1(x):=(\nabla b)_{ij}(x)=(\partial_jb^i(x))_{ij}$, and for $n\geq 2$,
\begin{eqnarray*}
B_n(x):=(b\cdot \nabla) B_{n-1}(x)-(\nabla b\cdot B_{n-1})(x).
\end{eqnarray*}

Now we can give our main result, which will be proven in Section 3.
\bt\label{Th2}
Assume that $(b, A)$ satisfy (\ref{EE1})-(\ref{EE3}) and H\"ormander's condition (\ref{eq}) at each point $x\in\mR^d$. Then 
for any $t>0$, the law of $X_t(x)$ is continuous in variable $x$ with respect to the total variation distance.
In particular, the semigroup $(P_t)_{t>0}$ has the strong Feller property, i.e., for any $t>0$ and $f\in\cB_b(\mR^d)$,
$$
x\mapsto \mE f(X_t(x))\mbox{ is continuous.}
$$
\et

\br
If $\mathrm{Rank}(A)=d$, then we can take $H(x):=|x|^2+1$ so that (\ref{EE1}) becomes 
$$
x\cdot b(x)\leq \kappa_1(|x|^2+1).
$$
In this case, the strong Feller property holds for SDE (\ref{SDE}) (cf. \cite{Zh1} \cite{Wa-Wa}).
\er

The topic about the smoothness of the distributional density of SDEs with jumps
has been studied for a long time since the work of Malliavin \cite{Ma}. We mention the following results:
\begin{itemize}
\item By using Girsanov's transformation,  Bismut  in \cite{Bi} established an integration by parts formula for Poisson functionals 
and then used it to study the smoothness of the distributional density of nondegenerate SDEs with jumps. His idea was systematically developed in 
the monograph \cite{Bi-Ja-Gr}.
\item In \cite{Pi}, Picard introduced a difference operator argument and derive a new criterion about the smoothness of the distributional density of Poisson functionals.
Moreover, the criterion is also used to SDEs with jumps. Recently, Ishikawa and Kunita in \cite{Is-Ku} extended Picard's result to Wiener-Poisson functional cases.
Moreover, Cass \cite{Ca} studied the SDEs driven by Browian motions and Poisson point processes under H\"ormander's conditions.
However, the result in \cite{Ca} does not cover the cases of (\ref{eq}) and $\alpha$-stable noises.
\item If $b(x)=Bx$, condition (\ref{eq}) is also called Kalman's condition. In this case,  Priola and Zabczyk \cite{Pr-Za}
proved the existence of smooth density for the corresponding Ornstein-Uhlenbeck process. In \cite{Zh2}, X. Zhang proved the existence
of density for SDE (\ref{SDE}) when $b$ is smooth Lipschitz continuous. In special degenerate cases, the smoothness of the density is also obtained
(cf. \cite{Zh2, Zh3}).
\end{itemize}


To the best of the authors' knowledge, Theorem \ref{Th2} is the first result about the regularization effect of L\'evy noises under {\it full} H\"ormander's conditions.
One motivation of our studies comes from the following stochastic oscillators studied in \cite{Ec-Ha, Re-Th, Ca0} etc.:
\begin{align}\label{6-1}
\begin{cases}
\dif z_i(t)=u_i(t)\dif t,& i=1,\cdots,d,\\
\dif u_i(t)=-\partial_{z_i}H(z(t),u(t))\dif t, & i=2,\cdots,d-1,\\ 
\dif u_i(t)=-[\partial_{z_i}H(z(t),u(t))+\gamma_iu_i(t)]\dif t+\sqrt{T_i}\dif W^i_{S^i_t}, & i=1, d,
\end{cases}
\end{align}
where $d\geq 3$, $\gamma_1,\gamma_d\in\mR$, $T_1,T_d>0$, and
$$
H(z,u):=\sum_{i=1}^d\left(\frac{1}{2}|u_i|^2+V(z_i)\right)+\sum_{i=1}^{d-1}U(z_{i+1}-z_{i}).
$$
The typical examples of $V$ and $U$ are
$$
V(z)=\frac{|z|^2}{2},\ \ U(z)=\frac{|z|^2}{2}+\frac{|z|^4}{4}.
$$
The Hamiltonian $H$ describes a chain of particles with nearest-neighbor interaction.  We have
\bp\label{Pr1}
Assume that $V, U\in C^\infty(\mR)$ are nonnegative and $\lim_{|z|\to\infty}V(z)=\infty$ so that $H$ is a Lyapunov function.
If $U$ is strictly convex, then (\ref{EE1}), (\ref{EE3}) and (\ref{eq}) hold.
\ep

This proposition will be proven in Section 4.

\section{Proof of Theorem \ref{Th1}}

Below, we fix a point $\lambda_0\in\Lambda$ and a neighbourhood $E_{\lambda_0}$ of $\lambda_0$. We divide the proof into three steps.
\\
\\
(1) Let GL$(d)\simeq\mR^d\times\mR^d$ be the set of all $d\times d$-matrix. Define
$$
K_n:=\Big\{A\in \mathrm{GL}(d): \|A\|\leq n,\ \ \det(A)\geq1/n\Big\}.
$$
Then $K_n$ is a compact subset of GL$(d)$. Let $\Phi_n\in C^\infty(\mR^d\times\mR^d)$ be a smooth function so that
$$
\Phi_n|_{K_n}=1,\ \ \Phi_n|_{K^c_{n+1}}=0,\ \ 0\leq \Phi_n\leq 1.
$$
For each $\lambda\in\Lambda$ and $n\in\mN$, let us define a finite measure $\mu_{\lambda,n}(\dif x)$ by
$$
\mu_{\lambda,n}(A):=\mE\Big[1_A(X_\lambda )\Phi_n(\Sigma^X_\lambda )\Big],\ \ A\in\sB(\mR^d).
$$
Then for each $\varphi\in C^\infty_b(\mR^d)$, by \cite[p.100, Proposition 2.1.4]{Nu}, we have
\begin{align*}
\int_{\mR^d}\nabla\varphi(x)\mu_{\lambda,n}(\dif x)&=\mE\Big[\nabla\varphi(X_\lambda )\Phi_n(\Sigma^X_\lambda )\Big]
=\mE\Big[\varphi(X_\lambda )\delta(\Phi_n(\Sigma^X_\lambda )(\Sigma^X_\lambda )^{-1}DX_\lambda )\Big],
\end{align*}
where $\nabla=(\p_1,\cdots,\p_d)$ and $\delta$ is the dual operator of $D$ (also called divergence operator). From this, by {\bf (H1)} and H\"older's inequality, we derive that
$$
\left|\int_{\mR^d}\nabla\varphi(x)\mu_{\lambda,n}(\dif x)\right|\leq \|\varphi\|_\infty C(\lambda,n),
$$
where $C(\lambda,n)$ is locally bounded in $\lambda$. Hence, $\mu_{\lambda,n}$ is absolutely continuous with respect to the Lebesgue measure
(cf. \cite{Nu}), and in particular, the density $p_{\lambda,n}$ satisfies 
$$
\int_{\mR^d}|\nabla p_{\lambda,n}(x)|\dif x\leq C(\lambda,n),
$$
which implies that $p_{\lambda,n}$ is locally bounded in $\mW^{1,1}(\mR^d)$ with respect to $\lambda$.
By Rellich-Kondrachov's compact embedding theorem (cf. \cite[p.168, Theorem 6.3]{Ad}), $\{p_{\lambda,n}\}_{\lambda\in E_{\lambda_0}}$ is compact in $L^1_{loc}(\mR^d)$, and
by Fr\'echet-Kolmogorov's theorem (cf. \cite[Ch 10]{Yo}), we have
\begin{align}
\lim_{|y|\to 0}\sup_{\lambda\in E_{\lambda_0}}\int_{B_M}|p_{\lambda,n}(x)-p_{\lambda,n}(x+y)|\dif x=0,\label{MB0}
\end{align}
where $B_M:=\{x\in\mR^d: |x|\leq M\}$ and $M>0$.
\\
\\
(2) Let $\phi\in C^\infty_c(B_1)$ be a nonnegative smooth function with $\int\phi=1$. For $\eps>0$, let 
$$
\phi_\eps(x):=\eps^{-d}\phi(\eps^{-1}x).
$$
For $f\in \cB_b(\mR^d)$ with support in $B_M$, let
$$
f_\eps(x):=\int_{\mR^d}f(y)\phi_\eps(x-y)\dif y.
$$
Noticing that
\begin{align*}
&\mE[(f(X_\lambda)-f_\eps(X_{\lambda}))\Phi_n(\Sigma^X_{\lambda})]
=\int_{\mR^d}(f(y)-f_\eps(y))p_{\lambda,n}(y)\dif y\\
&\qquad=\int_{\mR^d}f(y)\int_{\mR^d}(p_{\lambda,n}(y)-p_{\lambda,n}(y-x))\phi_\eps(x)\dif x\dif y,
\end{align*}
and in view of $f|_{B^c_M}=0$, we have
\begin{align}
|\mE[(f(X_\lambda)-f_\eps(X_{\lambda}))\Phi_n(\Sigma^X_{\lambda})]|
&\leq\|f\|_\infty\int_{B_M}\int_{\mR^d}|p_{\lambda,n}(y)-p_{\lambda,n}(y-x)|\phi_\eps(x)\dif x\dif y\no\\
&\leq\|f\|_\infty\sup_{x\in B_\eps}\int_{B_M}|p_{\lambda,n}(y)-p_{\lambda,n}(y-x)|\dif y.\label{MB1}
\end{align}
On the other hand, since $\mD^{k,q}=(I-\cL)^{-k}(L^q)$ by Meyer's inequality for any $q>1$, 
where $\cL=-\delta D$ is the Ornstein-Uhlenbeck operator, by the interpolation inequality, we have
$$
\|DX_\lambda-DX_{\lambda_0}\|_q\leq C\|X_\lambda-X_{\lambda_0}\|^{\frac{1}{2}}_q\|X_\lambda-X_{\lambda_0}\|^{\frac{1}{2}}_{2,q},
$$
which together with {\bf (H1)} and {\bf (H2)} implies that for any $q\in(1,p)$,
$$
\lim_{\lambda\to\lambda_0}\|DX_\lambda-DX_{\lambda_0}\|_q=0.
$$
Hence,
\begin{align}
\mbox{$\lambda\to\Sigma^X_\lambda$ is continuous in probability.}\label{MB2}
\end{align}
Observe that
\begin{align*}
|\mE(f(X_\lambda)-f(X_{\lambda_0}))|&\leq|\mE(f(X_\lambda)-f_\eps(X_{\lambda}))|
+|\mE(f(X_{\lambda_0})-f_\eps(X_{\lambda_0}))|+\mE|f_\eps(X_\lambda)-f_\eps(X_{\lambda_0})|\\
&\leq|\mE[(f(X_\lambda)-f_\eps(X_{\lambda}))\Phi_n(\Sigma^X_{\lambda})]|
+2\|f\|_\infty\mE|1-\Phi_n(\Sigma^X_{\lambda})|\\
&\quad+|\mE[(f(X_{\lambda_0})-f_\eps(X_{\lambda_0}))\Phi_n(\Sigma^X_{\lambda_0})]|
+2\|f\|_\infty\mE|1-\Phi_n(\Sigma^X_{\lambda_0})|\\
&\quad+\|f\|_\infty\int_{B_M}\mE|\phi_\eps(X_\lambda-y)-\phi_\eps(X_{\lambda_0}-y)|\dif y.
\end{align*}
By (\ref{MB0}), (\ref{MB1}),  (\ref{MB2}) and taking limits in order $\lambda\to\lambda_0$, $\eps\to 0$ and $n\to\infty$, we obtain
\begin{align}
\lim_{\lambda\to\lambda_0}\sup_{\|f\|_\infty\leq 1, f|_{B^c_M}=0}|\mE(f(X_\lambda)-f(X_{\lambda_0}))|
\leq 4\lim_{n\to\infty}\mP(\Sigma^X_{\lambda_0}\notin K_n)\stackrel{\bf (H3)}{=}0.\label{MB3}
\end{align}
(3) Lastly, noticing that for any $M>0$,
\begin{align*}
\sup_{\|f\|_\infty\leq 1}|\mE(f(X_\lambda)-f(X_{\lambda_0}))|
&\leq\sup_{\|f\|_\infty\leq 1, f|_{B^c_M}=0}|\mE(f(X_\lambda)-f(X_{\lambda_0}))|\\
&\quad+\mP(|X_\lambda|>M)+\mP(|X_{\lambda_0}|>M),
\end{align*}
by (\ref{MB3}), Chebyshev's inequality and {\bf (H1)}, we get
$$
\lim_{\lambda\to\lambda_0}\sup_{\|f\|_\infty\leq 1}|\mE(f(X_\lambda)-f(X_{\lambda_0}))|=0.
$$
The proof is thus completed by {\bf (H1)}, {\bf (H3)} and \cite[p.92, Theorem 2.1.1]{Nu}.
\section{Proof of Theorem \ref{Th2}}

The following lemma is proven in \cite[Lemma 2.1]{Zh2}.
\bl\label{EU1}
For $s>0$, set $\Delta\ell^j_s:=\ell^j_s-\ell^j_{s-}$ and
\begin{align*}
\mS_0:=\{\ell\in\mS: \{s: \Delta\ell^j_s>0\} \mbox{ is dense in $[0,\infty)$},\forall j=1,\cdots,m\}.
\end{align*}
Under (\ref{ER}), we have $\mu_\mS(\mS_0)=1$.
\el
Fix $\ell\in\mS_0$ and consider the following SDE:
\begin{align}
\dif X^\ell_t(x)=b(X^\ell_t(x))\dif t+A\dif W_{\ell_t},\ \ X^\ell_0=x.\label{SDE1}
\end{align}
The following result is proven in \cite[Theorem 3.1]{Zh3}.
\bt
Under (\ref{EE1})-(\ref{EE3}), there exists a unique solution to SDE (\ref{SDE1}) so that for all $t>0$,
\begin{align}
\mE\left[\exp\left\{\frac{2\sup_{s\in[0,t]}H(X^\ell_s(x))}{\e^{\kappa_1 t}(\kappa_2|\ell_t|+1)}\right\}\right]
\leq C_{\kappa_2,\kappa_3}\e^{H(x)},\label{EE93}
\end{align}
where $C_{\kappa_2,\kappa_3}\geq 1$. In particular, we have
$$
\mE f(X_t(x))=\mE(\mE f(X^\ell_t(x))|_{\ell=S}).
$$
\et
For proving the conclusion of Theorem \ref{Th2}, by Lemma \ref{EU1}, it suffices to show that for each $\ell\in\mS_0$ and $t>0$,
\begin{align}
\mbox{the law of $X^\ell_t(x)$ is continuous in $x$ with respect to the total variation norm.}\label{EY2}
\end{align}
For any $n\in\mN$, let $\chi_n(x)$ be a cut-off function on $[0,\infty)$ with
$$
\chi_n|_{B_n}=1,\ \ \chi_n|_{B^c_{n+1}}=0,\ \ 0\leq\chi_n\leq 1,
$$
and set
$$
b_n(x)=b(x)\chi_n(H(x)).
$$
Since $H\in C^\infty(\mR^d;\mR_+)$ and $\lim_{|x|\to\infty}H(x)=\infty$, we have
$$
b_n\in C^\infty_b(\mR^d).
$$
Consider the following SDE:
\begin{align}
\dif X^n_t(x)=b_n(X^n_t(x))\dif t+A\dif W_{\ell_t},\ \ X^n_0=x.\label{SDE2}
\end{align}
For fixed $t>0$ and $n\in\mN$, it is easy to see that {\bf (H1)} and {\bf (H2)} hold for $x\mapsto X^n_t(x)$.
On the other hand, the Malliavin covariance matrix of $X^n_t(x)$ has the following expression (cf. \cite[Lemma 4.5]{Zh3}):
$$
\Sigma^{X^n_t}_x=J^n_t(x)\left(\sum_{k=1}^m\int^t_0K^n_s(x) a_{\cdot k}
(K^n_s(x)a_{\cdot k})^*\dif \ell^k_s\right)(J^n_t(x))^*,
$$
where $J^n_t(x)$ and $K^n_t(x)$ solve the following matrix valued ODE:
$$
J^n_t(x)=I+\int^t_0\nabla b_n(X^n_s(x))\cdot J^n_s(x)\dif s
$$
and
$$
K^n_t(x)=I-\int^t_0J^n_s(x)\cdot \nabla b_n(X^n_s(x))\dif s.
$$
Define
$$
B^H_n:=\Big\{x\in\mR^d: H(x)<n\Big\}.
$$
If $(b,A)$ satisfies H\"ormander's condition (\ref{eq}) at point $x\in B^H_n$, then it is easy to see that
$(b_n,A)$ also satisfies H\"ormander's condition (\ref{eq}) at point $x\in B^H_n$.
Thus, from the proof of \cite[Theorem 1.1]{Zh2}, one sees that $\Sigma^{X^n_t}_x$ is invertible almost surely for $x\in B^H_n$. 
Using Theorem \ref{Th1}, for any $y\in B^H_n$, we have
\begin{align}
\lim_{x\to y}\sup_{\|f\|_\infty\leq 1}|\mE [f(X_t^n(x))-f(X_t^n(y))]|=0.\label{EY1}
\end{align}
 
Now, for any $x\in B^H_n$, define a stopping time
$$
\tau_n(x):=\inf\Big\{t\geq 0: H(X^\ell_t(x))\geq n\Big\}.
$$
 By the uniqueness of the solution to SDE, we have
$$
X^n_t(x)=X^\ell_t(x), \forall t<\tau_n(x),\ a.s.
$$
Let $f$ be a bounded nonnegative measurable function. For any $x,y\in B^H_n$, we have
\begin{align*}
|\mE[f(X^\ell_t(x))-f(X^\ell_t(y))]|&\leq|\mE[f(X^\ell_t(x))1_{t<\tau_n(x)}-f(X^\ell_t(y))1_{t<\tau_n(y)}]|\\
&\quad+\|f\|_\infty\mP(t\geq \tau_n(x))+\|f\|_\infty\mP(t\geq \tau_n(y))\\
&=|\mE[f(X^n_t(x))1_{t<\tau_n(x)}-f(X^n_t(y))1_{t<\tau_n(y)}]|\\
&\quad+\|f\|_\infty\mP(t\geq \tau_n(x))+\|f\|_\infty\mP(t\geq \tau_n(y))\\
&\leq|\mE[f(X^n_t(x))-f(X^n_t(y))]|\\
&\quad+2\|f\|_\infty\mP(t\geq \tau_n(x))+2\|f\|_\infty\mP(t\geq \tau_n(y)).
\end{align*}
Hence, by (\ref{EY1}) and (\ref{EE93}), we obtain
\begin{align*}
&\lim_{x\to y}\sup_{\|f\|_\infty\leq 1}|\mE[f(X^\ell_t(x))-f(X^\ell_t(y))]|\leq 4\lim_{n\to\infty}\sup_{|x-y|\leq 1}\mP(t\geq \tau_n(x))\\
&\qquad\leq 4\lim_{n\to\infty}\sup_{|x-y|\leq 1}\mP\left(\sup_{s\in[0,t]}H(X^\ell_s(x))\geq n\right)\\
&\qquad\leq 4\lim_{n\to\infty}\frac{1}{n}\sup_{|x-y|\leq 1}\mE\left(\sup_{s\in[0,t]}H(X^\ell_s(x))\right)=0.
\end{align*}
The proof is complete.
\section{Proof of Proposition \ref{Pr1}}

Let $x=(z_1,\cdots,z_d,u_1,\cdots,u_d)\in\mR^d\times\mR^d$ and define
$$
b(x):=b(z,u):=  \Big(u_1,\cdots, u_d, -[\partial_{z_1}H+\gamma_1u_1],\cdots,-\partial_{z_i}H,\cdots,-[\partial_{z_d}H+\gamma_du_d]\Big)
$$
and
$$
{\mbox{$A=(a_{i,j})$ with $a_{d+1,d+1}=\sqrt{T_1}$, $a_{2d,2d}=\sqrt{T_d}$, $a_{i,j}=0$ for other $i,j$.}}
$$
Clearly,
$$
b(x)\cdot\nabla H (x)=-\gamma_1^2u_1^2-\gamma^2_du_d^2\leq 0.
$$
Moreover,
$$
\sum_i\p_i H(x)a_{i,d+1}=\sqrt{T_1}u_1,\ \ \sum_i\p_i H(x)a_{i,2d}=\sqrt{T_d}u_d
$$
and
$$
\sum_{ij}\p_i\p_j H(x)a_{i,d+1}a_{j,d+1}=T_1,\ \ \sum_{ij}\p_i\p_j H(x)a_{i,2d}a_{j,2d}=T_d.
$$
Hence, (\ref{EE1}) and (\ref{EE3}) hold.

Let us now check (\ref{eq}). Let $\sV(x)$ be a vector field defined by
\begin{align*}
\sV(x):=\sV(z,u)&:=\sum_{i=1}^db_i(z,u)\p_{z_i}+\sum_{i=1}^db_{i+d}(z,u)\p_{u_i}\\
&= \sum_{i=1}^d u_i\partial_{z_i}-\big(\gamma_1 u_1+  V'(z_1)-U'(z_2-z_1)\big)\partial_{u_1}\\
&-\sum_{i=2}^{d-1}\big( V'(z_i)-U'(z_{i+1}-z_i)+U'(z_i-z_{i-1})\big) \partial_{u_i}\\
 &\quad- \big(\gamma_d u_d+ V'(z_d)+U'(z_d-z_{d-1}) \big)\partial_{u_d}.
\end{align*}
Here the prime denotes the differential.
Set $\sU_0:=\partial_{u_1}$ and define recursively
$$
\sU_n:=[\sU_{n-1},\sV]=\sU_{n-1}\sV-\sV\sU_{n-1},\ \ n\in\mN.
$$
By direct calculations, we have
$$
\sU_1=\p_{z_1}-\gamma_1\p_{u_1},
$$
$$
\sU_2= U''(z_2-z_1)\partial _{u_2}+\big(\gamma_1^2 -V''(z_1)-U''(z_2-z_1)\big)\partial_{u_1}-\gamma_1 \partial_{z_1}
$$
and
\begin{align*}
  \sU_3&=U''(z_2-z_1)\partial _{z_2}+\big(\gamma_1^2 -V''(z_1)-U''(z_2-z_1)\big)\partial_{z_1}
\\ &\quad+\left(\gamma_1V''(z_1)+\gamma_1 U''(z_2-z_1)+u_1V^{(3)}(z_1)+(u_2-u_1) U^{(3)}(z_2-z_1)\right)\partial_{u_1}
\\ &\quad+ \left((u_1-u_2) U^{(3)}(z_2-z_1)-\gamma_1U''(z_2-z_1)\right)\partial_{u_2}.
\end{align*}
By the induction, it is easy to see that for any $k=1,\cdots,d-2$,
$$
\left\{
\begin{aligned}
\sU_{2k}&=U''(z_{k+1}-z_{k})\cdots U''(z_2-z_1)\p_{u_{k+1}}+\sum_{i=1}^k(f_{ki}(x)\p_{z_i}+g_{ki}(x)\p_{u_i}),\\
\sU_{2k+1}&=U''(z_{k+1}-z_{k})\cdots U''(z_2-z_1)\p_{z_{k+1}}+\sum_{i=1}^k(\tilde f_{ki}(x)\p_{z_i}+\tilde g_{ki}(x)\p_{u_i})+h_k(x)\p_{u_{k+1}},
\end{aligned}
\right.
$$
where $f_{ki}, g_{ki},\tilde f_{ki},\tilde g_{ki}, h_k$ are smooth functions.
Since $U''>0$, we have
\begin{align}
\p_{u_1},\p_{z_1},\cdots,\p_{u_{d-1}},\p_{z_{d-1}}\in\mathrm{Span}\{\sU_0,\sU_1,\cdots,\sU_{2d-3}\}.\label{ERT1}
\end{align}
On the other hand, since
$$
[\p_{u_d}, \sV]=\p_{z_d}-\gamma_d\p_{u_d},
$$
by (\ref{ERT1}) we further have
$$
\p_{u_1},\p_{z_1},\cdots,\p_{u_{d}},\p_{z_{d}}\in\mathrm{Span}\Big\{\sU_0,\sU_1,\cdots,\sU_{2d-3},\p_{u_d}, [\p_{u_d},\sV]\Big\},
$$
which means that (\ref{eq}) holds.

\vspace{5mm}

{\bf Acknowledgements:}

The authors deeply thank the referees for their very useful suggestions so that Theorem 1.1 are improved.
Zhao Dong and Xuhui Peng are supported   by 973 Program, No. 2011CB808000 and
Key Laboratory of Random Complex Structures and Data Science, No. 2008DP173182, NSFC, Nos.: 10721101, 11271356, 11371041.
Xicheng Zhang is supported by NNSFs of China (Nos. 11271294, 11325105).

\end{document}